\documentclass[12pt,notitlepage,twoside,a4paper]{amsart}
 \usepackage{amsfonts}

\usepackage{amsmath,amssymb,enumerate}

\usepackage{epsfig,fancyhdr,color}

\usepackage{amssymb}
\usepackage{amsmath,amsthm}
\usepackage{latexsym}
\usepackage{amscd}
\usepackage{psfrag}
\usepackage{graphicx}
\usepackage[latin1]{inputenc}
\usepackage[all]{xy}
\usepackage[mathcal]{eucal}

\definecolor{NoteColor}{rgb}{1,0,0}

\renewcommand{\textsc}{\textcolor{red}}

\newtheorem*{theorem 1}{\rm\bf Proposition 1}
\newtheorem*{theorem 2}{\rm\bf Proposition 2}

\theoremstyle{definition}

\theoremstyle{remark}

\def\interieur#1{\mathord{\mathop{\kern 0pt #1}\limits^\circ}}

\title[A Commentary on Teichm\"uller's paper]{A Commentary on Teichm\"uller's paper \emph{Bestimmung der extremalen quasikonformen Abbildungen bei geschlossenen orientierten Riemannschen Fl\"achen}}

\author[A. A'Campo-Neuen, N. A'Campo, V. Alberge, A. Papadopoulos]{Annette A'Campo-Neuen, Norbert A'Campo,
\\ Vincent Alberge and Athanase Papadopoulos}

\thanks{The final version of this paper will appear in the \emph{Handbook of Teichm\"uler theory}, Vol. V. The authors were partially supported by the French ANR program FINSLER. V. Alberge and A. Papadopoulos acknowledge support from the U.S. National Foundation grants DMS 1107452, 1107263, 1107367 ``RNMS  GEometric structures And Representation varieties." V. Alberge and A. Papadopoulos  would like to thank 
the Graduate Center and Hunter College of the City University of New York, where part of this work was done.}

 \address{
  Universit\"at Basel,  Mathematisches Institut,\\
Spiegelgasse, No. 1, CH-4051 Basel, Switzerland\\
Institut de Recherche Math\'ematique Avanc\'ee\\ CNRS et Universit\'e de Strasbourg\\\small 7 rue Ren\'e
  Descartes - 67084 Strasbourg Cedex, France\\
  The Graduate Center, City University of New York,\\
365 Fifth Avenue,
New York, NY 10016 USA. \\
  \tt{annette.acampo@unibas.ch\\ Norbert.ACampo@unibas.ch}\\
    \tt{alberge@math.unistra.fr\\
   papadop@math.unistra.fr 
}}

\date{\today}

\begin{document}

  \maketitle
  \begin{abstract}
  
This is a mathematical commentary on Teichm\"uller's paper 
  \emph{Bestimmung der extremalen quasikonformen Abbildungen bei geschlossenen orientierten Riemannschen Fl\"achen} (Determination of extremal quasiconformal maps of closed oriented Riemann surfaces) \cite{T29}, (1943). This paper is among the last (and may be the last one) that Teichm\"uller wrote on the theory of moduli. It contains the proof of the so-called Teichm\"uller existence\index{extremal quasiconformal mappings!existence}\index{theorem!Teichm\"uller existence} theorem for a closed surface of genus $g\geq 2$. For this proof, the author defines a mapping between a space of equivalence classes of marked Riemann surfaces (the Teichm\"uller space) and a space of equivalence classes of certain Fuchsian groups (the so-called Fricke\index{Fricke space} space).\footnote{The consideration of such a space of Fuchsian groups, with the definition of a topology on that space, is already contained in the papers of Poincar\'e and Klein, cf. \cite{acampo&ji&papadop} and the references there to the papers by Poincar\'e and Klein. The Fricke space\index{Fricke space} has acquired in the recent literature  a new status, as a subspace of the character variety of the representation space of the fundamental group of the surface in the Lie group $\mathrm{PSL}(2,\mathbb{R})$.} After that, he defines a map between the latter and the Euclidean space of dimension $6g-6$. Using Brouwer's theorem of invariance of domain, he shows that this map is a homeomorphism. This involves in particular a careful definition of the topologies of Fricke space, the computation of its dimension, and comparison results between hyperbolic distance and quasiconformal dilatation. The use of the invariance of domain theorem is in the spirit of Poincar\'e and Klein's use of the so-called ``continuity principle"\index{continuity method}  in their attempts to prove the uniformization theorem.
\end{abstract}

  \tableofcontents
\section{Introduction}\label{intro}

We comment on the paper \cite{T29} by Oswald Teichm\"uller. In his article, the author proves a result which Ahlfors called the ``existence theorem\index{extremal quasiconformal mappings!existence} for extremal quasiconformal mappings,"  which is also known as ``Teichm\"uller's existence theorem" for arbitrary closed orientable surfaces of genus $g\geq 2$, and a related result which says that the Teichm\"uller space of a closed surface of genus $g$ is homeomorphic to $\mathbb{R}^{6g-6}$. (In fact, the two theorems are consequences of each other, but in Teichm\"uller's paper, the former is obtained as a consequence of the latter.) The article is a sequel to  \cite{T20} in which this result was presented as a ``conjecture," valid for a general surface of finite type. In the present paper, the author says that he restricts to the case of closed surfaces ``in order to make the essential points clear," and that ``the more general statement shall be proved in a later article." 

The paper \cite{T29} was published in 1943, and it is among the last ones that Teichm\"uller wrote. At the end of this paper, the author refers to the results of the paper \cite{T32}, which appeared in 1944 (the year after Teichm\"uller's death) \footnote{Teichm\"uller was drafted in the German army in July 1939, and never came back to university. He was killed at the Eastern Front in 1943. From 1939 until his death, he continued working on mathematics, but he could write only during his spare time. In his paper \cite{T23} (see also the commentary \cite{T23C}),  Teichm\"uller writes : ``Because I only have a limited vacation time at my disposal, I cannot give reasons for many things, but only assert."} in the ephemeral journal \textit{Deutsche Mathematik}\footnote{This journal was founded by L. Bieberbach, and it was published between 1936 and 1944.} together with his three other papers \cite{T30}, \cite{T31} and \cite{T33}.\footnote{Presumably, the last published paper by Teichm\"uller is \cite{T34}.}

Concerning the article which is our subject here, let us quote Ahlfors, from the comments he made in his \emph{Collected Works} (vol. II p. 1), concerning his paper \cite{Ahlfors}.  He writes about Teichm\"uller's existence theorem:\index{extremal quasiconformal mappings!existence}  ``I found this proof rather hard to read''  but ``I did not doubt its validity.'' He also explains, concerning his own  proof of that theorem: ``My attempted proof on these lines had a flaw, and even my subsequent correction does not convince me today.'' He refers to Bers, and he says that the latter gave  in \cite{B1}  ``a very clear version of Teichm\"uller's proof.''

In \cite{B1} Bers writes : ``Our arrangement of the arguments preserves the logical structures of Teichm\"uller's proof; the details are carried out differently. More precisely, we work with the most general definition of quasiconformality, we rely on the theory of partial differential equations in some crucial parts of the argument, and we make use of a simple set of moduli for marked Riemann surfaces." In fact, in the paper \cite{B1}, written 17 years after the one of Teichm\"uller, Bers uses results which a priori were not known to Teichm\"uller, in particular a result of Morrey \cite{Morrey}, which is a weak version of the so-called \emph{Ahlfors-Bers} theorem,\index{theorem!Ahlfors-Bers}\index{theorem!measurable Riemann mapping} also known as the \textit{measurable Riemann mapping theorem}. However, in another paper, \cite{B2}, written 16 years after \cite{B1} (and which is among the very last papers that Bers wrote),\footnote{This paper appears as the antepenultimate in Bers' \emph{Selected works} edition.} Bers insists on the fact that all the arguments given by Teichm\"uller in his paper \cite{T29} are correct. The paper \cite{B2} is considered by Bers as a postscript to his paper  \cite{B1}. In the abstract of \cite{B2}, Bers says that he gives a ``simplified version of Teichm\"uller's proof (independent of the theory of Beltrami equations with measurable coefficients)\footnote{Bers means Morrey's theory.} of a proposition underlying his continuity argument\index{continuity method}  for the existence\index{extremal quasiconformal mappings!existence} part of his theorem on extremal quasiconformal mappings." He then writes in the introduction: 
\begin{quote}\small
In proving a basic continuity assertion (Lemma 1 in \S\,14C of \cite{B1}) I made use of a property of quasiconformal mappings (stated for the first time in \cite{B0} and also in \S\,4F of \cite{B1}) which belongs to the theory of quasiconformal mappings with bounded measurable Beltrami coefficients (and seems not to have been known to Teichm\"uller). Some readers concluded that the use of that theory was indispensable for the proof of Teichm\"uller's theorem. This is not so, and Teichm\"uller's own argument is correct. This argument can be further simplified and this simplified argument will be presented below. Then we will briefly describe Teichm\"uller's argument.
\end{quote}

Concerning the general presentation of \cite{T29},  Teichm\"uller states precisely his results, with relatively clear and detailed proofs (unlike some of his other papers, e.g.  \cite{T20} and \cite{T32}, and the commentaries \cite{T20C} and \cite{T32C}). However, Ahlfors declares in \cite{Ahlfors}, concerning Teichm\"uller's paper \cite{T29}, that the proof is an ``anticlimax'', in comparison with the same author's paper \cite{T20}, which is a ``brilliant and unconventional paper.'' 
The presentation here is very close to the one in \cite{T24}. This is not surprising because in the two papers, Teichm\"uller uses a continuity argument and applies Brouwer's theorem of invariance of domain.

It is probably time now to state Teichm\"uller's existence\index{extremal quasiconformal mappings!existence}\index{theorem!Teichm\"uller existence} theorem, which is the main result of the paper \cite{T29}. The theorem says that in any  homotopy (or isotopy) class of homeomorphisms between two closed Riemann surfaces, there exists a quasiconformal homeomorphism which realizes the infimum of the quasiconformal dilatation in the given homotopy class. The uniqueness of such a homeomorphism and  a precise description of it in terms of quadratic differentials were already proved in \cite{T20}.

Let us highlight some of the main ideas contained in the paper \cite{T29}, besides the proof of this existence theorem:\index{extremal quasiconformal mappings!existence}

\begin{enumerate}
\item The description of the topologies of various spaces of equivalence classes of Fuchsian groups, and in particular of a space which is called today Fricke space,\index{Fricke space} and the computation of their dimension.
\item The fact that Teichm\"uller space is homeomorphic to the Fricke space and that the two spaces are homeomorphic to $\mathbb{R}^{6g-6}$. 
\item A correct use of Brouwer's theorem of invariance of domain\index{theorem!invariance of domain} in the setting of Riemann surfaces.\footnote{The idea of this method goes back to Poincar\'e and Klein, who called this the ``method of continuity."\index{continuity method} However, these authors, in using this method, bumped over several difficulties, because they were applying it to Riemann's moduli space, which is not a manifold. Furthermore, a rigorous definition of the notion of dimension was missing. It was Brouwer who proved the needed result with all the needed rigour. For a survey see \cite{acampo&ji&papadop}.}
\item Estimates on the transformation of lengths of curves under a $K$-quasicon\-formal mapping between hyperbolic surfaces. These estimates are used in the proof of the continuity of the map that is involved in Brouwer's theorem of invariance of domain.
\end{enumerate}
 
The structure of Teichm\"uller's paper is very clear, with the following sections:

\begin{itemize}
\item ``Introduction," where the results are stated, the main result being a proof of the Teichm\"uller existence theorem.\index{extremal quasiconformal mappings!existence}
\item ``Topological determination and uniformization of the surface $\mathfrak{M}$," where the author recalls the notion of \emph{canonical dissection}\index{canonical dissection!surface} of a Riemann surface.\footnote{This notion was already considered by Poincar\'e; cf. the historical report \cite{acampo&ji&papadop}.} Using this notion and the uniformization theorem, he establishes a one-to-one correspondence between the space of equivalence classes of marked Riemann surfaces and a space of conjugacy classes of Fuchsian groups uniformizing a closed Riemann surface.
\item ``$\mathfrak{G}$ as a linear group," where several spaces of equivalence classes of Fuchsian groups are introduced and shown to be topological manifolds.
\item ``The extremal quasiconformal mappings," where the notion of extremal quasiconformal mapping (in modern terms, \emph{Teichm\"uller mapping}) and the relation with quadratic differentials are recalled from the paper \cite{T20}. Using these notions, Teichm\"uller reformulates the main result of the paper \cite{T20}. 
\item ``The fundamental continuity proof,"\index{continuity method}  where several fundamental results that involve at the same time quasiconformal geometry and hyperbolic geometry are obtained. The results concern the continuity of a map between the Fricke space\index{Fricke space} and a Euclidean space. These  continuity results are essential for the main result of the paper, namely, Teichm\"uller's existence theorem.\index{extremal quasiconformal mappings!existence}
\item ``Remark," where a fundamental inequality on the effect of hyperbolic length under a quasiconformal mapping is proved. The result complements the continuity result of the previous section. It is a consequence of the proofs in that section. The result shows the continuity of the map from Euclidean space to Fricke space at the origin, and it is also used in the proof that Fricke space is arcwise connected.
\item ``Construction of a sufficiently regular mapping from $\underline{\mathfrak{M}}$ to $\underline{\mathfrak{M}'}$," which concerns mappings between marked Riemann surfaces. The author proves that in each homotopy class of homeomorphism between two Riemann surfaces there exists a quasiconformal map.\footnote{Teichm\"uller  assumes, in this papers and others, that the quasiconformal homeomorphisms considered are ``sufficiently regular." In the present construction, they are real-analytic except at a finite set of arcs which are the sides of a triangulation parametrized by real-analytic arcs.}
\item ``The continuity proof,"\index{continuity method}  which contains the conclusion of the fact that the space denoted by $\mathfrak{E}$ is homeomorphic to $\mathbb{R}^{6g-6}$, based on the so-called ``continuity method," or, Brouwer's theorem of invariance of domain.\index{theorem!invariance of domain}
\item ``Connection to the theory of moduli," which is a concluding section, where the author makes the connection with the moduli problem as stated by Riemann.\index{problem!Riemann moduli} The author proves that $\underline{\mathfrak{R}}$ is homeomorphic to $\mathbb{R}^{6g-6}$, and he mentions relations with his other works on the subject and directions for further investigations.
\end{itemize}
Throughout this commentary, we shall keep the logical trend of Teichm\"uller's paper and  state some of the important results with the author's original wording.

 \section{Teichm\"uller's paper}

The introduction of the paper is brief and explicit. Teichm\"uller starts by saying that he is ``glad now to be able to actually prove'' the existence\index{extremal quasiconformal mappings!existence} of an extremal quasiconformal mapping for a given problem (which will be recalled later). This result had the status of a conjecture in the previous paper \cite{T20}. In the present paper, the author declares that the theorem will only be proved for closed orientable surfaces of genus $g>1$ and, unlike what he did in the paper \cite{T20}, he will use several times the uniformization theorem. He also recalls that the case $g=1$ was completely treated in \cite{T20} (\S\,25 to \S\,29). Teichm\"uller  declares that the proof will be based on a ``continuity argument.''\index{continuity method}  By this, he means that he will exhibit a mapping between $\mathbb{R}^{6g-6}$ and another space (in fact, the space which is called today \textit{Teichm\"uller space}, denoted here by $\underline{\mathfrak{R}}$, as in the paper \cite{T32})\footnote{In the paper \cite{T20}, this paper is denoted by $R^{\sigma}$, with $\sigma=6g-6$.} and he will show that this mapping is injective and continuous. He will then appeal to Brouwer's \textit{invariance of domain} theorem to conclude that this mapping is a homeomorphism. Teichm\"uller  says that he considers only closed orientable surfaces in order to avoid a technical topological theorem that concerns the homotopy equivalence relation between maps between surfaces. He says that to prove that theorem in the more general case would take an unproportional amount of space and distract the reader from the main point of the paper, which is the  ``continuity argument.'' In fact, the technical part would be to show, in the general case, that there exists a ``sufficiently regular mapping" in each homotopy class. This is harder for surfaces with distinguished points. He nevertheless promises to treat in the near future the general case, that is, the case of surfaces which are orientable or not, with or without boundary and with or without distinguished points. Let us also note that the proof of the existence theorem in the case of the torus that we mentioned above ($g=1$) does not use the method of continuity, but that Teichm\"uller gave a proof of his existence theorem, using the method of continuity,\index{continuity method}  in the case where the surface is a pentagon, that is, a disc with five distinguished points on the boundary, cf. \cite{T24}. In the last section of the paper  \cite{T29}, he says: ``During my research, I have always kept in mind the aim to give a continuity proof for the existence of quasiconformal mappings similar to the one in the case of the pentagon."

 The aim of the section called ``Topological determination and uniformization of the surface $\mathfrak{M}$'' is to give several equivalent definitions of the object that will be denoted 
 by  $\underline{\mathfrak{R}}$ in the rest of the paper. The author fixes a base surface  $\mathfrak{M}_0$ and defines an equivalence relation on the set of pairs $\left( \mathfrak{M}, H \right)$ where $\mathfrak{M}$ is a Riemann surface which has the same topological type as $\mathfrak{M}_0$ and $H : \mathfrak{M}_0 \rightarrow \mathfrak{M}$ an orientation-preserving homeomorphism. The equivalence relation is such that $\left( \mathfrak{M}, H \right)$ is equivalent to $\left( \mathfrak{M}^\prime, H^\prime \right)$ if $\mathfrak{M}=\mathfrak{M}^\prime$ and $H^\prime \circ H^{-1}$ is homotopic to the identity.\footnote{Teichm\"uller uses \textit{isotopy}, and he notes that by a result of Mangler (cf. \cite{Mangler}),  isotopy is equivalent to homotopy in the setting he considers.} He calls such an equivalence class a \textit{topologically determined Riemann surface}.\index{topologically determined!Riemann surface}\index{Riemann surface!topologically determined}\footnote{Interpreting the equality $\mathfrak{M}=\mathfrak{M}^\prime$ as meaning that $\mathfrak{M}$ is conformally equivalent to $\mathfrak{M}^\prime$  by a map $f:\mathfrak{M}\to\mathfrak{M}^\prime$ and the relation ``$H^\prime \circ H^{-1}$ is homotopic to the identity" as saying that the map $f \circ H$ is homotopic to $H'$, we get the  definition of Teichm\"uller space to which we are used today.} Teichm\"uller  denotes by $\underline{\mathfrak{M}}$ the equivalence class of pairs $\left( \mathfrak{M}, H \right)$. He then introduces the natural notions of topological  and of conformal maps between two topologically determined Riemann surfaces. Using this language, the Teichm\"uller space of $\mathfrak{M}_0$ is then the set $\underline{\mathfrak{M}}$ of ``topologically determined Riemann surfaces up to conformal maps.'' 
 
 Teichm\"uller then recalls the notion of \emph{canonical dissection}\index{canonical dissection!surface} of the surface $\mathfrak{M}_0$. This is a collection of $2g$ simple loops on this surface, $A_1,\ldots, A_g,B_1,\ldots B_g$, with the same base point and disjoint except at this base point, and satisfying, as elements of the fundamental group of the surface, the relation 
\begin{equation}\label{eq1}
\Pi_{i=1}^{g}{A_i B_i A_i ^{-1}B_i ^{-1}}=1,
\end{equation} 
where $1$ denotes the identity element of the group.
  The term \emph{dissection} is used because when the surface is cut along the union of these curves, we get a polygon with $4g$ sides. In particular, the set of homotopy classes of  loops $A_1,\ldots,A_g,B_1,\ldots B_g$ is a set of generators for the fundamental group of the surface $\mathfrak{M}_0$. A dissection\index{canonical dissection!surface} of $\mathfrak{M}_0$ induces a dissection of any topologically determined Riemann surface $(\mathfrak{M},H)$, via the topological map $H$,  which is well defined up to an inner automorphism of the fundamental group of $\mathfrak{M}$. Conversely, a dissection on $\mathfrak{M}$ defines a topological determination $H$ of that surface, that is, a topologically marked surface $(\mathfrak{M},H)$. Teichm\"uller appeals to a theorem of Mangler, contained in the same paper \cite{Mangler}, for the fact that a topologically determined surface defines a surface with a marking of its fundamental group, up to an inner automorphism.\footnote{This is related to the fact that the group of homotopy classes of homeomorphisms (that is, the mapping class goup) of the closed surface is in natural one-to-one correspondence with the outer automorphism group of the fundamental group, a result attributed to Dehn and Nielsen.}

Using the uniformization theorem, Teichm\"uller identifies the universal cover of the Riemann surface $\mathfrak{M}$ with upper half-plane $\mathbb{H}$ and the fundamental group of $\mathfrak{M}$ to a discrete subgroup of $\mathrm{PSL}(2,\mathbb{R})$ of orientation-preserving isometries of $\mathbb{H}$. With this, he obtains the following (we use his own words, in the translation of \cite{T29}):\begin{quote}
\textit{We have a one-to-one correspondence between, on the one hand, the topologically determined Riemann surfaces $\underline{\mathfrak{M}}$ that are given only up to conformal mappings and, on the other hand, certain collections of linear mappings $A_1 , \cdots , A_g , B_1 , \cdots , B_g$ of $\mathbb{H}$ onto itself that are determined only up to a common transformation with a linear mapping of $\mathbb{H}$ onto itself.}
\end{quote}
The term ``certain collections" refers to a property satisified by the elements $\left\lbrace A_i , B_i \right\rbrace_{1\leq i\leq g}$, namely, that they are different from the identity and are of hyperbolic type, that they satisfy the relation (\ref{eq1}) when the elements are considered as sitting in $\mathrm{PSL}(2,\mathbb{R})$, that they generate a discrete subgroup (denoted by $\mathfrak{G}$) which acts properly discontinuously on $\mathbb{H}$ and such that the quotient is a closed surface of genus $g$. Such a collection of elements is called an ``admissible collection"  and the set of all admissible collections forms a space which will be denoted by $\mathfrak{E}$. Teichm\"uller quotes another result of Mangler \cite{Mangler} which implies that the natural map between the space of topologically determined surfaces up to conformal maps (that is, Teichm\"uller space) and Fricke space is surjective.
 
%
%
%
%
%

The third section, titled ``$\mathfrak{G}$ as a linear group'', is about topology. Teichm\"uller introduces several spaces of equivalence classes of matrices and he equips them with topologies, showing that these spaces are topological manifolds and computing their dimension. This will be used in the ``method of continuity"\index{continuity method}, or ``invariance of domain,"\index{theorem!invariance of domain} which is a crucial ingredient in the proof of the main result of the paper. 

Teichm\"uller introduces a space $\mathfrak{E}$ equipped with a topology induced by its inclusion in a $6g-6$-dimensional topological manifold $\mathfrak{D}$. 
This space $\mathfrak{E}$ is called today the \textit{Fricke space}.\index{Fricke space} An important fact which is shown in this section is that for every $\mathfrak{t}$ in $\mathfrak{E}$, there exists an element $\mathfrak{G}$ in $\mathfrak{C}$ generated by a collection of elements $A_1,A_2,B_1,B_2,\ldots,A_g,A_g,B_g,B_g$ such that a finite number of invariants $\lambda$ of these elements and of certain of their products determine $\mathfrak{t}$ uniquely and in a continuous manner. Let us see this in some detail. 

 Teichm\"uller recalls that to every hyperbolic element $A$ of $\mathrm{PSL}(2,\mathbb{R})$ one can associate an attractive fixed point, a repulsive fixed point and an invariant, which is
 the real number $\lambda>1$ such that $A$ is conjugate to the transformation $z\mapsto \lambda  z$ of the upper half-plane.\footnote{The value of $\lambda$ is the exponential of the hyperbolic length of the simple closed curve associated to $A$.} He briefly recalls the notion of an $n$-dimensional topological manifold and he introduces three sets, 
 $\mathfrak{A}$, $\mathfrak{B}$ and $\mathfrak{C}$. 
 
 The set $\mathfrak{A}$ is the set of all $6g$-tuples of elements of  $\mathrm{PSL}(2,\mathbb{R})$. This is a $6g$-dimensional manifold. The set $\mathfrak{B}$ is the subset of $6g$-tuples of $\mathfrak{A}$ that satisfy the relation (\ref{eq1}) and such that the first two elements, $A_1$ and $B_1$ are hyperbolic with four distinct fixed points. The set $\mathfrak{C}$ is the set of admissible collections $A_1 , \cdots , A_g , B_1 , \cdots , B_g$ of elements $\mathrm{PSL}(2,\mathbb{R})$.
 
  Applying the \textit{implicit function theorem} to 
  \[D=\Pi_{i=1}^{g}{A_i B_i A_i ^{-1}B_i ^{-1}},
  \] seen as a map $\mathfrak{A}\rightarrow\mathrm{PSL}(2,\mathbb{R})$, Teichm\"uller shows that $\mathfrak{B}$ is a $6g-3$-dimensional manifold. He then defines the spaces $\mathfrak{D}$ and $\mathfrak{E}$ which are the quotients of $\mathfrak{B}$ and $\mathfrak{C}$ respectively by the natural actions on them of the group $\mathrm{PSL}(2,\mathbb{R})$ by conjugation.  
  
  He concludes with the following:
  \begin{quote}
  \textit{The set $\underline{\mathfrak{R}}$ of all classes of conformally equivalent topologically determined Riemann surfaces of genus $g$ is in one-to-one correspondence with $\mathfrak{E}$.}
  \end{quote}
  The first set that is considered here is what became known later on as Teichm\"uller space\index{Teichm\"uller space} and the second space  is Fricke space.\index{Fricke space} Let us note that it is at this stage that  Teichm\"uller introduces the notation $\underline{\mathfrak{R}}$.  
  
Teichm\"uller equips $\mathfrak{D}$ with the topology induced by $\mathfrak{B}$  and he shows that this topology is Hausdorff, after introducing several technical lemmas. He  then shows that this space is a  $(6g-6)$-topological manifold and he announces the result saying that $\mathfrak{E}$ is an open connected subset of $\mathfrak{D}$ which is homeomorphic to $\mathbb{R}^{6g-6}$. He says that this statement ``was not known previously.'' He ends this section by showing that the elements of $\mathfrak{E}$ depend continuously on a finite number of invariants. Let us state his result in his own terms:
\begin{quote}
\textit{There is a finite set of a priori listable expressions in the $A_i$, $B_i$, as e.g. $A_1$, $A_1 B_1$, $\left( A_1  B_1\right)^2 B_2$, etc. with the following property: The class $\mathfrak{k}$ of the collection $\mathfrak{a}$ is uniquely determined by the invariants $\lambda_k$ of those finitely many hyperbolic mappings, and for every neighborhood $\mathfrak{U}$ of $\mathfrak{k}$ there is a $\delta>0$ such that the classes $\mathfrak{k}^\prime$ of all collections $\mathfrak{a}^\prime$ whose invariants $\lambda_k ^\prime$ satisfy the inequalities $\left( \lambda'_k-\lambda_k \right)<\delta$, lie in the neighborhood $\mathfrak{U}$.}
\end{quote}
Here, $\mathfrak{k}\in\mathfrak{E}$ and $\mathfrak{a}\in\mathfrak{C}$.\footnote{In his proof, Teichm\"uller uses a representative of $\mathfrak{k}$ in $\mathfrak{C}$ which is sometimes called \textit{normalized Fuchsian model} (cf. \cite{Imayoshi} p. 47). Let us also note that, in his proof, Teichm\"uller shows that ``$\mathfrak{k}$ is uniquely and continuously determined by'' precisely $24g-20$ invariants, whereas in \cite{Imayoshi} (Proposition 6.17), there is a similar   result with only $14g-11$ invariants. It seems obvious to the authors of the present report that Teichm\"uller was aiming directly to the proof of the main theorem, and did not try to obtain the best estimates for the intermediate results.}

The next section is called ``The extremal quasiconformal mappings." The author recalls a construction he presented in detail in \cite{T20} of what is called today the \textit{Teichm\"uller maps}, i.e. the maps referred to in the title of this section. Such a map is determined by a holomorphic quadratic differential (which the author calls an ``\textit{everywhere finite quadratic differential}'')\footnote{The word ``finite" originates in the fact that the area of the surface, for the singular flat metric associated to the quadratic differential, is finite. This is a condition on the poles. It is true that in the present paper of Teichm\"uller, all quadratic differentials considered are holomorphic and therefore the area is always finite, but the terminology originates in his paper \cite{T20}, and in that paper, the quadratic differentials are allowed to have poles at distinguished points.} $d\zeta^2$ on  $\mathfrak{M}_0$ and a real number $K\geq 1$. We denote such a map $E\left( K, d\zeta^2 \right)$.\footnote{Teichm\"uller uses a slightly different notation, by taking first a basis for the vector space of quadratic differentials.} Away from the zeros of $d\zeta^2$, a Teichm\"uller map is affine with respect to \textit{natural coordinates} of $d\zeta^2$, with constant ``dilatation quotient''\footnote{In the literature, this quantity is usually called \textit{maximal dilatation}.} equal to $K$. Teichm\"uller also recalls the fact that such a map is extremal and unique in its homotopy class. We note by the way that unlike the existence theorem, the proofs of the uniqueness theorem that were given in the first three decades after Teichm\"uller's work are modeled on those that he gave in (\cite{T20} Chapter 27).\footnote{We can quote here Ahlfors from \cite{Ahlfors}, talking about his own paper: ``[...] a complete proof of the uniqueness part of Teichm\"uller's theorem was included. Like all other known proofs of the uniqueness it was modeled on Teichm\"uller's own proof [...].''} Such a map naturally defines a new element of $\underline{\mathfrak{R}}$ denoted by $\underline{\mathfrak{M}}\left( K, d\zeta^2 \right)$. This element depends on $d\zeta^2$ ``up to a common \textit{positive} factor." In other words, a positive constant multiple of $d\zeta^2$, together with the same $K$, define the same extremal quasiconformal map. Using the bijective correspondence that he obtained between $\underline{\mathfrak{R}}$ and $\mathfrak{E}$, Teichm\"uller defines a mapping between $\mathbb{R}^{6g-6}$ and $\mathfrak{E}$. The element of $\mathfrak{E}$ which corresponds to $\underline{\mathfrak{M}}\left( K, d\zeta^2 \right)$ is denoted by $\mathfrak{e}\left( K, d\zeta^2 \right)$. Teichm\"uller concludes this section by announcing the following key result:
\begin{quote}
$\mathfrak{e}\left( K, d\zeta^2 \right)$ \textit{depends continuously on} $(K, d\zeta^2)$.
\end{quote}
He stresses the fact that continuity makes sense because a topology has been defined on $\mathfrak{D}$ and therefore on $\mathfrak{E}$. 

In the next section, called\index{continuity method}  ``The fundamental continuity proof,'' Teichm\"uller proves the result we just mentioned. Even though the details are technical, the idea is rather clear. Indeed, according to the section ``$\mathfrak{G}$ as a linear group'', we only need to show that the invariants corresponding to the group elements that determine the point $\mathfrak{e}\left( K, d\zeta^2 \right)$ depend continuously on $K$ and $d\zeta^2$. For this, the author considers two points of $\underline{\mathfrak{R}}$; $\underline{\mathfrak{M}}\left( K, d\zeta^2  \right)$ and $\underline{\mathfrak{M}}\left(  \tilde{K}, d\tilde{\zeta}^2 \right)$,  determined respectively by the Teichm\"uller maps $E\left( K, d\zeta^2 \right)$ and $E\left( \tilde{K}, d\tilde{\zeta}^2 \right)$ such that $\tilde{K}$ is close to $K$ and $d\tilde{\zeta}^2$ is close to $d\zeta^2$. He then lifts the composition of the two Teichm\"uller maps $E\left( \tilde{K}, d\tilde{\zeta}^2\right) E\left( K, d\zeta^2 \right)^{-1}$ to a quasiconformal mapping of $\mathbb{H}$ which conjugates the normalized Fuchsian group $\mathfrak{G}$ to a normalized Fuchsian group $\mathfrak{\tilde{G}}$. Let us recall that these two groups determine respectively $\underline{\mathfrak{M}}\left( K, d\zeta^2 \right)$ and $\underline{\mathfrak{M}}\left( \tilde{K}, d\tilde{\zeta}^2 \right)$. He concludes using a notion of \textit{average of the local dilatation}\footnote{This notion is contained in Teichm\"uller's paper, and the name is given to it by Bers in \cite{B2}.} that the invariants that characterize $\underline{\mathfrak{M}}\left( \tilde{K}, d\tilde{\zeta}^2 \right)$ are close to the invariants characterizing $\underline{\mathfrak{M}}\left( K, d\zeta^2 \right)$. It is interesting to observe that Teichm\"uller uses, for the needs of the average of the local dilatation, an idea already present in the paper \cite{T20}, \S\,35, which permits, among other things, to compare hyperbolic length with quasiconformal dilatation and which leads to the so-called Wolpert inequality (cf. \cite{wolpert}). 

In the section titled ``Remark,'' Teichm\"uller recalls the idea of the proof of the result of the preceding section and he states a consequence of this proof which will be useful in particular to show that $\mathfrak{E}$ is arcwise connected. This result will also be useful to show the continuity at the origin of the map between $\mathbb{R}^{6g-6}$ and the Fricke space, which we already alluded to. Let us state this idea explicitly.
\begin{quote}
\textit{Let $\underline{\mathfrak{M}}^\prime$ be a topologically determined Riemann surface of genus $g>1$, corresponding to the point $\mathfrak{e}^\prime \in\mathfrak{E}$. Then for every neighborhood $\mathfrak{U}$ of $\mathfrak{e}^\prime$ in $\mathfrak{D}$ there exists a $\delta>0$ with the following property: If a quasiconformal mapping $A$ from $\underline{\mathfrak{M}}^\prime$ onto a second topologically determined surface $\underline{\tilde{\mathfrak{M}}}$ has a dilatation quotient $D$ satisfying $D\leq 1+\delta$ everywhere, then the point $\tilde{\mathfrak{e}}$ of $\mathfrak{E}$ belonging to $\tilde{\mathfrak{M}}$ lies in $\mathfrak{U}$.}
\end{quote}
Let us only note that to prove this fact it is not necessary to consider what we called the average of the local dilatation.

The next section is called ``Construction of a sufficiently regular mapping from $\underline{\mathfrak{M}}$ to $\underline{\mathfrak{M}}^\prime$." Let us explain the term  ``sufficiently regular." Teichm\"uller starts by defining the notion of a ``\textit{sufficiently regular triangulation}'' of a Riemann surface $\mathfrak{M}$, as a triangulation whose vertices are parametrized by maps which are real analytic except at the  vertices. Then, a map between $\underline{\mathfrak{M}}$ and $\underline{\mathfrak{M}}^\prime$ is said to be sufficiently regular if it is a homeomorphism which transforms a sufficiently regular triangulation of $\underline{\mathfrak{M}}$ into a sufficiently regular triangulation of $\underline{\mathfrak{M}}^\prime$ and which is real analytic in the interior of the triangles. Furthermore, he imposes that the dilatation quotient is bounded at the vertices.\footnote{Regarding the existence of a quasiconformal map between two points of $\underline{\mathfrak{R}}$, we can quote Bers,  from his paper \cite{B1}: ``Had we demanded that the homeomorphism $f$ be continuously differentiable everywhere, the proof would be somewhat laborious. Since we use a very general definition of quasiconformality, the proof [that there exists a quasiconformal homeomorphism between two arbitrary marked Riemann surfaces] presents no difficulties and may be omitted." We note that Teichm\"uller considered maps which are differentiable except may be on a finite set of arcs, and that Bers' remark holds also in this more general case.}
 
In the section called ``The continuity proof,''\index{continuity method}  Teichm\"uller recalls and describes in detail the  continuous and injective map between $\mathbb{R}^{6g-6}$ and $\mathfrak{E}$. He denotes the image by $\mathfrak{E}^*$. Given that this set is included in a  $6g-6$-dimensional manifold, he uses Brouwer's theorem of invariance of domain\index{theorem!invariance of domain} in order to conclude that the map is a homemorphism. This proves in particular that the set $\mathfrak{E}^*$ is open in $\mathfrak{D}$. Then, using the existence of a sufficiently regular mapping between two points of $\mathfrak{E}$, using the notion of Riemannian metrics on surfaces representing conformal structures which he developed in his paper \cite{T20}, \S\,16, he shows that $\mathfrak{E}$ is arcwise connected. Finally, using the arcwise connectedness of $\mathfrak{E}$, he concludes this section by showing that $\mathfrak{E}^* =\mathfrak{E}$ and therefore that this space is homeomorphic to $\mathbb{R}^{6g-6}$. 
 
In the last section, called ``Connection to the theory of modules,'' Teichm\"uller establishes relations between the present work and his previous works. He also recalls his motivation, namely, the solution of the so-called \emph{problem of moduli}.\index{problem!Riemann moduli}\index{Riemann moduli problem} This problem originates in Riemann's work, and it consists in giving a precise meaning to Riemann's statement that the set of equivalence classes of Riemann surfaces, where the equivalence relation is that of conformal equivalence, has $3g-3$ complex moduli. More generally, the problem became that of understanding the space of moduli.  We refer the reader to the paper  \cite{acampo&ji&papadop} for the history of this problem. Teichm\"uller considers that the  \emph{problem of moduli} consists in turning the space $\underline{\mathfrak{R}}$ into an analytic manifold. He recalls that he outlined a solution to that problem in his paper \cite{T32} (see also the commentary in \cite{T32C}).

Teichm\"uller, in this last section, states several results of which he does not give the details, and which make connections between the work done in this paper and his other papers.  It is interesting to note that in the paper  \cite{T32}, in which Teichm\"uller presents a version of Teichm\"uller space equipped with a complex-analytic structure, he declares that he is not sure that this space (which he also denoted by $\underline{\mathfrak{R}}$)\footnote{In \cite{T32}, Teichm\"uller space is denoted by $\underline{\mathfrak{R}}$, and, as a space equipped with its complex analytic structure, it is denoted by $\underline{\mathfrak{C}}$.} coincides with the Teichm\"uller space that he studied in his paper \cite{T20} and  whose study is continued in the paper \cite{T29} that we are commenting on here. He writes, in particular, in  \cite{T32}: ``The space $\underline{\mathfrak{R}}$ consists of at most countably many connected parts. I believe that  $\underline{\mathfrak{R}}$ is in fact simply connected." In the paper \cite{T29}, there is no more questioning about this fact, and Teichm\"uller knows that the spaces considered are the same. In particular, the space $\underline{\mathfrak{R}}$  is connected and is homeomorphic to a $(6g-6)$-dimensional Euclidean space.

 Furthermore, Teichm\"uller says that the bijective correspondence between $\underline{\mathfrak{R}}$ and $\mathfrak{E}$ that was given at the beginning of the paper is continuous, when we consider on $\underline{\mathfrak{R}}$ the topology induced by the \textit{Teichm\"uller metric}. This completes the statement known today as the  \textit{Teichm\"uller theorem}, namely, that $\underline{\mathfrak{R}}$ is homeomorphic to $\mathbb{R}^{6g-6}$. 
 
\section{Some other approaches}
We mention briefly some later (but still old) approaches to Teichm\"uller's existence theorem. We already mentioned the works of Ahlfors \cite{Ahlfors}  and Bers \cite{B1}, \cite{B2}. Gerstenhaber and  Rauch inaugurated an approach using the theory of minimal surfaces, more precisely, the method of minimizing the so-called energy integral, or Douglas-Dirichlet functional \cite{GR1} \cite{GR2}. Hamilton in \cite{Ham} gave a proof of the existence theorem based on the geometry of Banach spaces (the spaces of integrable Beltrami differentials and the space of all Beltrami differentials) and maps between them. The proof works for general Riemann surfaces with or without distinguished points. The Gerstenhaber-Rauch approach is surveyed in \cite{DW}, together with other approaches. Krushkal' also developed a variational proof, cf. \cite{K1} and the book \cite{K2}. There are several modern proofs and generalizations of Teichm\"uller's existence theorem, in particular to infinite-dimensional Teichm\"uller spaces, and it is not possible to mention them here.

\end{document}